\documentclass[12pt]{article}
\catcode`\@=11
\def\demo{\noindent{\bf Proof .-}}
\def\section{\@startsection {section}{1}{\z@}{-3.5ex plus -1ex
minus-.2ex}{2.3ex plus .2ex}{\normalsize\bf}}

\newtheorem{lemma}{Lemma}
\newtheorem{proposition}{Proposition}

\newtheorem{example}{Example}
\newtheorem{remark}{Remark}

\newtheorem{corollary}{Corollary}
\def\lcm{{\rm lcm}\,}

\begin{document}
\begin{center}
{\LARGE\bf{A note on monomial ideals}}
\end{center}
\vskip.5truecm
\begin{center}
{Margherita Barile\\ Dipartimento di Matematica, Universit\`{a} di Bari, Via E. Orabona 4,\\70125 Bari, Italy\\e-mail:barile@dm.uniba.it}
\end{center}
\vskip1truecm
\noindent
{\bf Abstract} We show that the number of elements generating a squarefree monomial ideal up to radical can always be bounded above in terms of the number of its minimal monomial generators and the maximal height of its minimal primes. \vskip0.5 truecm
\noindent

\section*{Introduction}
Let $R$ be a commutative Noetherian ring with non-zero identity. If $I$ is an ideal of $R$, we say that the elements 
$a_1,\dots, a_s\in R$ {\it generate} $I$ {\it up to radical} if ${\rm Rad}\,(a_1,\dots, a_s)={\rm Rad}\,I$. The minimum number $s$ with this property is called the {\it arithmetical rank} of $I$.
If $\mu(I)$ is the minimum number of generators for $I$, then obviously
\begin{equation}\label{ara}{\rm ara}\,I\leq\mu(I).\end{equation}
\noindent
This inequality yields a trivial upper bound for the arithmetical rank of $I$.  Moreover, it is well-known that
\begin{equation}\label{height}h(I)\leq {\rm ara}\,I,\end{equation}
\noindent where $h(I)$ is the height of $I$. If equality holds in (\ref{ara}) and in (\ref{height}), $I$ is called a {\it complete intersection}. 
The following result by Schmitt and Vogel \cite{SV}, p. 249, shows that (\ref{ara}) can be improved if special combinatorial conditions on the generating set are fulfilled.
\begin{lemma}\label{Schmitt} Let $P$ be a finite subset of elements of $R$. Let $P_0,\dots, P_r$ be subsets of $P$ such that
\begin{list}{}{}
\item[(i)] $\bigcup_{i=0}^rP_i=P$;
\item[(ii)] $P_0$ has exactly one element;
\item[(iii)] if $p$ and $p'$ are different elements of $P_i$ $(0<i\leq r)$ there is an integer $i'$ with $0\leq i'<i$ and an element in $P_{i'}$ which divides $pp'$.
\end{list}
\noindent
We set $q_i=\sum_{p\in P_i}p^{e(p)}$, where $e(p)\geq1$ are arbitrary integers. We will write $(P)$ for the ideal of $R$ generated by the elements of $P$.  Then we get
$${\rm Rad}\,(P)={\rm Rad}\,(q_0,\dots,q_r).$$
\end{lemma}
We will apply  Lemma \ref{Schmitt} in the polynomial ring $R=K[x_1,\dots, x_n]$, where $K$ is a field, in order to determine upper bounds for the arithmetical rank of ideals generated by monomials, i.e., by products of indeterminates. Since ideals with the same radical have the same arithmetical rank, without loss of generality we can restrict our investigation to radical monomial ideals. These are the so-called {\it squarefree} monomial ideals, i.e., the ideals that are generated by  products of pairwise distinct indeterminates. \newline
 The aim of this paper is to replace inequality (1) by a more precise one in the case where $I$ is a squarefree monomial ideal. We will derive a closed expression, depending only on $\mu(I)$ and the heights of the minimal primes of $I$, which can be put on the right-hand of inequality (\ref{ara}) in order to obtain a better general upper bound.  This new upper bound is still sharp when $I$ is a complete intersection. Moreover it improves the general upper bound conjectured by Lyubeznik in \cite{L2} when $I$ has a small number of generators. \newline
Recall that every monomial ideal $I$ has a unique monomial generating set of cardinality $\mu(I)$: it is the one formed by the monomials in the ideal that are minimal with respect to the divisibility relation.
\section{A general upper bound}
Let $I$ be a squarefree monomial ideal of $R$, and let $M$ be the set of $\mu(I)$ monomial generators for $I$. For all $i=1,\dots, n$ let
\begin{equation}\label{mi}M_i=\{f\in M\mid x_i\mbox{ divides }f\}.\end{equation}
\noindent
Let $I_1, \dots, I_r$ be the minimal primes of $I$, so that $I=\cap_{j=1}^rI_j$. For all $j=1,\dots, r$ set
\begin{equation}\label{nuj}\nu_j=\max\{\vert M_i\vert\mid x_i\in I_j\},\end{equation}
and let $x_{i_j}\in I_j$ be such that $\nu_j=\vert M_{i_j}\vert$. Here the bars indicate the cardinality. Finally, let
\begin{equation}\label{nu}\nu(I)=\min\{\nu_j\mid j=1,\dots, r\}.\end{equation}
\noindent From  definitions (\ref{mi}), (\ref{nuj}) and (\ref{nu}) it follows that $\nu(I)\leq\vert M\vert=\mu(I)$. 
\begin{proposition}\label{mylemma} ara\,$I\leq\mu(I)-\nu(I)+1$.
\end{proposition}
\demo For the sake of simplicity we put $\mu=\mu(I)$ and $\nu=\nu(I)$. We are going to explicitly construct a set of at most $\mu-\nu+1$ polynomials which generate $I$ up to radical. For all $k=1,\dots, \mu$ let
$$L_k=\{\prod_{i=1}^kf_i\mid f_1, \dots, f_k\mbox{ are pairwise distinct elements of $M$}\}.$$
\noindent
Further, let
$$p_0=\prod_{j=1}^rx_{i_j}\in\displaystyle \cap_{j=1}^rI_j=I.$$
\noindent
Then put
\begin{equation}\label{p0} P_0=\{p_0\},\end{equation}
\noindent and, for all $i=1,\dots, \mu-\nu$,
$$P_i=L_{\mu-\nu+1-i}.$$
\noindent We show that, for all $i=1,\dots, \mu-\nu$,
\begin{equation}\label{condition}\mbox{ the product of any two distinct elements of $P_i$ is divisible by some element of $P_{i-1}$}.
\end{equation}
\noindent
 We first show it for $i>1$. Set $k=\mu-\nu+1-i$ and let $p,p'$ be distinct elements of $P_i=L_k$. Then there are two different $k$-subsets $\{f_1,\dots, f_k\}$ and $\{f_1',\dots, f_k'\}$ of $M$ such that
$$p=\prod_{i=1}^kf_i\qquad\qquad\mbox{ and }\qquad\qquad p'=\prod_{i=1}^kf_i'.$$
\noindent Up to a change of indices we may assume that $f_1'\not\in\{f_1,\dots, f_k\}$. Then $pp'$ is divisible by 
$$\prod_{i=1}^kf_i\cdot f_1'\in L_{k+1}=P_{i-1}.$$
Now consider two distinct elements $p,p'$ of $P_1=L_{\mu-\nu}$. By the above argument, $pp'$ is divisible by a product $\pi$ of at least $\mu-\nu+1$ distinct elements of $M$. By definition, for all $j=1,\dots, r$, exactly $\nu_j$ elements of $M$ are divisible by $x_{i_j}$, in particular, by definition (\ref{nu}), at least $\nu$ elements of $M$ are divisible by $x_{i_j}$, or, equivalently, at most $\mu-\nu$ are not divisible by $x_{i_j}$.  It follows that the product $\pi$ necessarily involves, as a factor, an element of $M$  which is divisible by $x_{i_j}$. Hence $pp'$ is divisible by $l_0$. 

\noindent This proves condition (\ref{condition}) and implies that $P_0, \dots, P_{\mu-\nu}$ fulfill the assumption (iii) of Lemma \ref{Schmitt} with $i'=i-1$. Moreover, we have that $P_{\mu-\nu}=L_1=M$, and $P_i\subset (M)$ for all $i=0,\dots,\mu-\nu$, so that $\left(\cup_{i=0}^{\mu-\nu}P_i\right)=(M)=I$. Thus assumption (i) of Lemma \ref{Schmitt} is satisfied by $P=M$. Finally, assumption (ii) is trivially true by (\ref{p0}). Therefore, 
if we set
$$q_i=\sum_{p\in P_i}p\qquad\qquad\mbox{ for }i=0,\dots, \mu-\nu,$$
\noindent
by Lemma \ref{Schmitt} we have that 
$$I=(M)={\rm Rad}\,(q_0,\dots, q_{\mu-\nu}),$$
\noindent
which implies the claim.
\begin{remark}\label{remark1} {\rm The construction of the sets $P_i$ given in the proof of Lemma \ref{mylemma} does not give rise to an efficient algorithm, since it requires the computation of $2^{\mu(I)}-1$ polynomial products. The output can be simplified if one reduces each set $P_i$ (with $i\geq1$) to the set of its elements which are minimal with respect to the divisibility relation. This will not affect condition (\ref{condition}), and will possibly produce a smaller number of polynomials $q_i$.  Sets $P_i$ fulfilling (7) and assumptions (i) and (ii) of Lemma \ref{Schmitt} can also be constructed by the following method, which is taken from \cite{B}, and is more convenient from a computational point of view. We first put $\Gamma_1=M$, and for all indices $i\geq2$, we recursively define $\Gamma_i$ as the set of all elements of 
$$\{\lcm(f,g)\mid f, g\in \Gamma_{i-1}, f\ne g\}$$
\noindent
which are are minimal with respect to the divisibility relation.
The procedure stops at step $N$ as soon as there is some $p_0\in I$ which divides all elements of $\Gamma_N$. At this point we set $P_0=\{p_0\}$ and $P_i=\Gamma_{N-i}$  for all $i=1,\dots, N-1$. Then the above condition (\ref{condition}) is trivially true.\newline The different approaches discussed here are not equivalent: one can find examples where they yield distinct numbers of polynomials generating $I$ up to radical.} 
\end{remark}
\par\bigskip\noindent Next we show that number $\mu(I)$ can be compared with $\nu(I)$. Let 
$$\tau(I)=\max\{h(I_j)\mid j=1,\dots, r\}.$$
\noindent Then we have
\begin{proposition}\label{mylemma2} $\mu(I)\leq\nu(I)\tau(I)$.
\end{proposition}
\demo We adopt the simplified notation used in the preceding proof and we also set $\tau=\tau(I)$.
In view of (\ref{nu}) we have to prove that $\mu\leq\nu_j\tau$ for all $j=1,\dots, r$. We show it for $j=1$. Up to renaming the indeterminates we may assume that $I_1=(x_1,\dots, x_s)$ for some $s\leq\tau$. Since $M\subset I_1$, by (\ref{mi}) we have $M=\cup_{i=1}^s M_i$, whence
$$\mu=\vert M\vert \leq \sum_{i=1}^s\vert M_i\vert\leq s\nu_1\leq\tau\nu_1,$$
\noindent
where the second inequality follows from (\ref{nuj}) for $j=1$. This completes the proof.
\par\bigskip\noindent
The next result is an immediate consequence of Propositions \ref{mylemma} and \ref{mylemma2}.  
\begin{corollary}\label{mycorollary} Let $J$ be an ideal of $R$ such that $I=$\,Rad\,$(J)$ is a (squarefree) monomial ideal. Then \begin{equation}\label{ara2}{\rm ara}\,J\leq \mu(I)-\displaystyle\frac{\mu(I)}{\tau(I)}+1.\end{equation}
\end{corollary}
\begin{remark}\label{ci}{\rm  The new upper bound given in Corollary \ref{mycorollary} is, like the one in (\ref{ara}), sharp if $I=J$ and $I$ is a complete intersection: in that case $\mu(I)=h(I)$ and  it is well-known that $I$ is pure-dimensional, i.e., all its minimal primes have the same height (see, e.g., \cite{BH}, Corollary 5.1.5), so that $h(I)=\tau(I)$.
}\end{remark}
\begin{remark}\label{jaballah}{\rm For a squarefree monomial ideal $I$, the number $\mu(I)$ can be bounded above in terms of  the number $n$ of indeterminates of $R$ and/or of the heights of the minimal primes of $I$. Some of these formulas can be found in \cite{J} or in  \cite{J0}. If we replace them in (\ref{ara2}) we obtain similar upper bounds for ara\,$I$. For example, Lemma (4) (a) in \cite{J} yields
$${\rm ara}\,I\leq \left(1-\frac1{\tau(I)}\right){n\choose{\lbrack\frac{n}2\rbrack}}+1,$$
\noindent
where the square brackets denote the integer part, and, according to  Satz (5) (a) in \cite{J}, we have
$${\rm ara}\,I\leq \left(1-\frac1{\tau(I)}\right){n\choose{h(I)-1}}+1,$$
\noindent whenever $h(I)\geq\frac{n+1}2$.
}\end{remark}
\section{On a conjecture by Lyubeznik}
It was conjectured by Lyubeznik in \cite{L2} that for every pure-dimensional ideal $J$ of $R$ of height $t$, 
\begin{equation}\label{Lyubeznik}{\rm ara}\,J\leq n-\left[\frac{n-1}{t}\right],\end{equation}
so that, in particular, for every pure-dimensional monomial ideal $I$,
\begin{equation}\label{Lyubeznik2}{\rm ara}\,I\leq n-\left[\frac{n-1}{\tau(I)}\right].\end{equation}
\noindent Inequality (\ref{Lyubeznik}) has been proven by the same author in \cite{L3} for all ideals $J$ in the localized ring $R_{(x_1,\dots, x_n)}$ (if $K$ is infinite), whereas (\ref{Lyubeznik2}) was established in \cite{L1}, Theorem 6, for monomial ideals $I$ with $\tau(I)=2,3$ (see also \cite{B}, Section 2, for a different proof in the case $\tau(I)=2$).  Other special cases were examined in \cite{SeV} and \cite{SV}. In \cite{B} it is conjectured that (\ref{Lyubeznik2}) holds for every squarefree monomial ideal $I$.
The results in the previous section allow us to show, and even improve, inequality (\ref{Lyubeznik2}) for a new class of monomial ideals.  Corollary \ref{mycorollary} implies the following:
\begin{corollary}\label{mycorollary2} Let $J$ be an ideal of $R$ such that $I=$\,Rad\,$(J)$ is a monomial ideal generated by at most $n-1$ (squarefree) monomials. Then 
$${\rm ara}\,J\leq n-\frac{n-1}{\tau(I)}.$$ 
\end{corollary}
\begin{remark}{\rm In the class of squarefree monomial ideals $I$ such that $\mu(I)\leq n-1$   the upper bounds given in Proposition \ref{mylemma} and in Corollary \ref{mycorollary} are in general strictly better than the one in  (\ref{Lyubeznik2}), as the next example shows.
}\end{remark}
\begin{example}\label{example1}{\rm Suppose that char\,$K=0$ and in $R=K[x_1,\dots, x_9]$ consider the ideal
$$I=(x_1x_2,\ x_1x_3,\ x_2x_4,\ x_4x_5,\ x_4x_6,\ x_2x_7,\ x_6x_8,\ x_6x_9),$$
\noindent
with $\mu(I)=8$. Its minimal primes can  be quickly computed by CoCoA \cite{CoC}:
$$I_1=(x_1, x_2, x_4, x_6),\qquad I_2=(x_2, x_3, x_4, x_6),\qquad I_3=(x_1, x_2, x_5, x_6),$$
$$I_4=(x_2, x_3, x_5, x_6),\qquad I_5=(x_1, x_4, x_6, x_7),\qquad I_6=(x_1, x_2, x_4, x_8, x_9),$$
$$I_7=(x_2, x_3, x_4, x_8, x_9),\qquad I_8=(x_1, x_4, x_7, x_8, x_9).$$
\noindent Thus we have $\tau(I)=5$, so that $n-\lbrack\frac{n-1}{\tau(I)}\rbrack=9-\lbrack\frac85\rbrack=8$. Hence (\ref{Lyubeznik2}) yields the trivial inequality (\ref{ara}).  From Corollary \ref{mycorollary} we derive that ara\,$I\leq 8-\frac85+1=7.4$, i.e., ara\,$I\leq 7$. Moreover,
$$\vert M_3\vert=\vert M_5\vert=\vert M_7\vert=\vert M_8\vert=\vert M_9\vert=1,\qquad\vert M_1\vert =2,$$
$$\vert M_2\vert=\vert M_4\vert=\vert M_6\vert=3,$$
\noindent so that, by definition (\ref{nuj}), $\nu_j=3$ for all $j=1,\dots, 8$. Therefore, by definition (\ref{nu}), $\nu(I)=3$. Hence Proposition \ref{mylemma} gives us ara\,$I\leq 8-3+1=6$.   The method from \cite{B} that we mentioned in Remark \ref{remark1} also yields ara\,$I\leq6$. In fact ara\,$I=5$: on the one hand, by virtue of Lemma \ref{Schmitt}, $I$ is generated up to radical by the 5 polynomials
$$x_1x_2,\ x_1x_3+x_2x_4,\ x_2x_7+x_4x_6,\ x_4x_5+x_6x_8,\ x_6x_9;$$
\noindent
on the other hand, according to \cite{L4}, Theorem 1, and \cite{H}, Example 2, p.~414, we have pd\,$(R/I)\leq$\,ara\,$I$, where pd denotes the projective dimension, and CoCoA tells us that here pd\,$(R/I)=5$. }
\end{example} 
In the previous example,  the best upper bound for the arithmetical rank is the one derived from Proposition \ref{mylemma}. Sometimes this proposition yields the exact value of the arithmetical rank. An interesting example of this kind can be obtained by a slight modification of the ideal in Example \ref{example1}.  
\begin{example}{\rm Suppose that char\,$K=0$ and in $R=K[x_1,\dots, x_9]$ consider the ideal
$$I=(x_1x_2,\ x_1x_3,\ x_2x_4,\ x_4x_5,\ x_4x_6,\ x_6x_7,\ x_6x_8,\ x_6x_9),$$
\noindent
with $\mu(I)=8$. The minimal primes of $I$ are
$$I_1=(x_1, x_4, x_6),\qquad I_2=(x_2, x_3, x_4, x_6),\qquad I_3=(x_1, x_2, x_5, x_6),$$
$$I_4=(x_2, x_3, x_5, x_6),\qquad I_5=(x_1, x_4, x_7, x_8, x_9),\qquad I_6=(x_2, x_3, x_4, x_7, x_8, x_9).$$
\noindent Now $\tau(I)=6$, so that $n-\lbrack\frac{n-1}{\tau(I)}\rbrack=9-\lbrack\frac86\rbrack=8$. As in the previous example, (\ref{Lyubeznik2}) yields a trivial upper bound.  From Corollary \ref{mycorollary} we deduce that ara\,$I\leq 8-\frac86+1=7.\bar6$, i.e., ara\,$I\leq 7$.  Since  $\nu(I)=3$, by Proposition \ref{mylemma} we again have ara\,$I\leq6$. The same result is obtained if we use the method from \cite{B}, and also if we apply Lemma \ref{Schmitt} directly to the minimal monomial generators of $I$: $I$ is generated up to radical by the 6 polynomials
$$x_2x_4,\ x_1x_2+x_4x_6,\ x_4x_5+x_6x_7,\ x_1x_3,\ x_6x_8,\ x_6x_9.$$
\noindent
This time  pd\,$(R/I)=6$, so that  ara\,$I=6$. }
\end{example}

\end{document}